\documentclass[11pt,reqno]{amsart}

\usepackage{a4wide}
\usepackage{amsmath}
\usepackage{amsfonts}
\usepackage{amssymb}
\usepackage{graphicx}

\usepackage[utf8]{inputenc}
\usepackage{hyperref}
\usepackage{esint}
\usepackage{pdfsync}

\newcommand{\mean}{\mathbb E}
\newcommand{\prob}{\mathbb P}
\newcommand{\de}{{\rm d}}

\newtheorem{Theorem}{Theorem}
\newtheorem{Lemma}{Lemma}
\newtheorem{Proposition}{Proposition}

\theoremstyle{remark}
\newtheorem{Remark}{Remark}

\newcommand{\R}{\mathbb R}
\newcommand{\C}{\mathbb C}
\newcommand{\borel}{\mathcal B}
\newcommand{\real}{\mathcal Re\,}
\newcommand{\im}{\mathcal Im\,}

\newcommand{\beq}{\begin{equation}}
\newcommand{\eeq}{\end{equation}}
\newcommand{\blem}{\begin{Lemma}}
\newcommand{\elem}{\end{Lemma}}
\newcommand{\bthm}{\begin{Theorem}}
\newcommand{\ethm}{\end{Theorem}}
\newcommand{\bprop}{\begin{Proposition}}
\newcommand{\eprop}{\end{Proposition}}
\newcommand{\bcor}{\begin{Corollary}}
\newcommand{\ecor}{\end{Corollary}}
\newcommand{\brmk}{\begin{Remark}}
\newcommand{\ermk}{\end{Remark}}

\makeatletter
\newcommand*{\mint}[1]{%
  \mint@l{#1}{}%
}
\newcommand*{\mint@l}[2]{%
  \@ifnextchar\limits{%
    \mint@l{#1}%
  }{%
    \@ifnextchar\nolimits{%
      \mint@l{#1}%
    }{%
      \@ifnextchar\displaylimits{%
        \mint@l{#1}%
      }{%
        \mint@s{#2}{#1}%
      }%
    }%
  }%
}
\newcommand*{\mint@s}[2]{%
  \@ifnextchar_{%
    \mint@sub{#1}{#2}%
  }{%
    \@ifnextchar^{%
      \mint@sup{#1}{#2}%
    }{%
      \mint@{#1}{#2}{}{}%
    }%
  }%
}
\def\mint@sub#1#2_#3{%
  \@ifnextchar^{%
    \mint@sub@sup{#1}{#2}{#3}%
  }{%
    \mint@{#1}{#2}{#3}{}%
  }%
}
\def\mint@sup#1#2^#3{%
  \@ifnextchar_{%
    \mint@sub@sup{#1}{#2}{#3}%
  }{%
    \mint@{#1}{#2}{}{#3}%
  }%
}
\def\mint@sub@sup#1#2#3^#4{%
  \mint@{#1}{#2}{#3}{#4}%
}
\def\mint@sup@sub#1#2#3_#4{%
  \mint@{#1}{#2}{#4}{#3}%
}
\newcommand*{\mint@}[4]{%
  \mathop{}%
  \mkern-\thinmuskip
  \mathchoice{%
    \mint@@{#1}{#2}{#3}{#4}%
        \displaystyle\textstyle\scriptstyle
  }{%
    \mint@@{#1}{#2}{#3}{#4}%
        \textstyle\scriptstyle\scriptstyle
  }{%
    \mint@@{#1}{#2}{#3}{#4}%
        \scriptstyle\scriptscriptstyle\scriptscriptstyle
  }{%
    \mint@@{#1}{#2}{#3}{#4}%
        \scriptscriptstyle\scriptscriptstyle\scriptscriptstyle
  }%
  \mkern-\thinmuskip
  \int#1%
  \ifx\\#3\\\else_{#3}\fi
  \ifx\\#4\\\else^{#4}\fi
}
\newcommand*{\mint@@}[7]{%
  \begingroup
    \sbox0{$#5\int\m@th$}%
    \sbox2{$#5\int_{}\m@th$}%
    \dimen2=\wd0 %
    \let\mint@limits=#1\relax
    \ifx\mint@limits\relax
      \sbox4{$#5\int_{\kern1sp}^{\kern1sp}\m@th$}%
      \ifdim\wd4>\wd2 %
        \let\mint@limits=\nolimits
      \else
        \let\mint@limits=\limits
      \fi
    \fi
    \ifx\mint@limits\displaylimits
      \ifx#5\displaystyle
        \let\mint@limits=\limits
      \fi
    \fi
    \ifx\mint@limits\limits
      \sbox0{$#7#3\m@th$}%
      \sbox2{$#7#4\m@th$}%
      \ifdim\wd0>\dimen2 %
        \dimen2=\wd0 %
      \fi
      \ifdim\wd2>\dimen2 %
        \dimen2=\wd2 %
      \fi
    \fi
    \rlap{%
      $#5%
        \vcenter{%
          \hbox to\dimen2{%
            \hss
            $#6{#2}\m@th$%
            \hss
          }%
        }%
      $%
    }%
  \endgroup
}

\begin{document}

\title[An extension of Hewitt's inversion formula]{An extension of Hewitt's inversion formula and its application to fluctuation theory}
\maketitle

\bigskip
\centerline{ \c Serban  E. B\u adil\u a\footnote{Supported by Project 613.001.017 of the Netherlands Organisation for Scientific Research (NWO)} \let\thefootnote\relax\footnote{\:\emph{Department of Mathematics and Computer Science, Eindhoven University of Technology, The Netherlands.} \\ \indent\:\emph{E-mail: serban.e.badila@gmail.com}}}

\bigskip
{\setlength{\parindent}{2pt}


\begin{abstract}
We analyze fluctuations of random walks with generally distributed increments. Integral representations for key performance measures are obtained by extending an inversion theorem of Hewitt \cite{Hewitt} for Laplace-Stieltjes transforms. Another important step in the analysis involves introducing the so-called harmonic measures associated to the walk.
 It is also pointed out that such representations can be explicitly calculated, if one assumes a form of rational structure for the increment transform. Applications include, but are not restricted to, queueing and insurance risk problems.
\end{abstract}

\vspace{.3cm}



\phantom{2010 Mathematics Subject Classification:~}Secondary {30E20}
}

\section{Introduction}

A standard assumption in the theory of fluctuations of random walks as they appear in, e.g., queueing or insurance applications, is that the increment of the random walk can be represented as the difference of two \emph{independent} random variables.
 In the context of queueing theory, random walks with increments which {\em do not have} this property appear as embedded at arrival epochs of customers in a single server queue in which the service requirement of the current customer is correlated with the time until the next arrival. In risk reserve processes that appear in insurance, the independence assumption is violated when the current claim size depends on the time elapsed since the previous arrival, and hence on the premium gained meanwhile. 

The purpose of the present paper is to show how much can be done for random walks which do not satisfy this independence assumption, regarding their maxima (as in the waiting time/maximum aggregate loss), their minima (idle periods/deficit at ruin) and their excursions which are related to busy periods or the time to ruin. 

From a queueing perspective, it turns out that the busy period is  a more sensitive issue to study than the idle period or the waiting time, as it appears from the proof of Theorem \ref{Theorem integral representation} below. For this purpose, we present a generalization of an inversion formula for Fourier-Stieltjes transforms due to Hewitt~\cite{Hewitt}, which in turn is an extension of P. L\'evy's inversion formula. All these results are essentially variations on the Dirichlet integral for complex-valued functions of bounded variation.

From a stochastic point of view, the information contained by the increments of the random walk is sufficient to infer about the extreme statistics; and since successive increments are independent, the usual form of Hewitt's inversion formula is sufficient to obtain the integral representations; this is how it was used originally by Spitzer \cite{Spitzer57} to derive the Laplace-Stieltjes transform of the maxima of partial sums. He also related these transforms to the Wiener-Hopf problem (see also Cohen \cite{CohenSingleServer}, Ch. II.5, for the relation with the Wiener-Hopf equation as it appears in Probability Theory). For the  derivation of the length of an excursion above the starting level, there is more information needed, namely, that given by the partial sums $\sum_{i=1}^n B_i$ {\em  together} with the partial sums of the embedded random walk $S_n= \sum_i^n B_i - \sum_i^n A_i$. It is possible to derive the excursion lengths still using Hewitt's formula when the $A_i$'s are independent of the $B_i$'s, and this was carried out in Kingman \cite{Kingman62a}. It is shown here that if one extends Hewitt's inversion formula, a similar derivation is possible for the case when there is  arbitrary correlation inside the vectors $(A_i,B_i)$, which means the random walk $S_n$ can have generally distributed increments.

Hewitt's approach was to find a most general inversion identity for Laplace-Stieltjes transforms with a view on Harmonic Analysis; this is more than needed for our purposes.
Instead of trying to find a most general instance of inversion, we will focus on obtaining a sufficiently broad result to apply to random walks as they appear for example in the study of workload/insurance related problems. One can then hope that the result itself will find applications in other related areas of probability and statistics.

In the queueing literature, Conolly \cite{Conolly60} obtained (the transform of) the busy period together with the number of customers served in a queue with exponential inter-arrivals and independent Erlang distributed service times. 
Conolly's results from \cite{Conolly60}  were  extended to general independent inter-arrivals and service times in Finch \cite{Finch} and in Kingman \cite{Kingman62a}.

The time to ruin has been studied in the insurance literature by deriving recursion formulae, typically obtained by discretising the claim sizes. For example, Dickson and Waters \cite{DicksonWaters} present various approximation methods for its numerical computation. 
Another good reference is Prabhu \cite{Prabhu61} $\S$3, who obtained an integral equation for the time to ruin starting from a positive capital $u$, in the Cramer-Lundberg risk reserve model; Borovkov and Dickson~\cite{BorovkovDikson} obtain series representations for the distribution of the time to ruin in the Sparre-Andersen risk reserve model with exponentially distributed claims and general renewal inter-arrivals. Besides these results, there exists a significant amount of literature on the Gerber-Shiu functions which contain the time to ruin as a special case. 

The paper is organized in the following way: In Section \ref{section inversion}, we extend Hewitt's inversion formula to allow for probability  distribution functions on $\R^2$ which do not have a product form (see Remark \ref{remark original Hewitt}). One of the ingredients of the proof consists of having a precise meaning for the conditional distribution of $B$ given $A$; this is settled as a preliminary.
The approach used for studying fluctuations of random walks involves obtaining integral representations for the above-mentioned quantities.
 The busy period, idle period, transient workload, and the related insurance functionals can still be determined in the form of a Cauchy integral, once Theorem \ref{generalized Hewitt} is combined with a version of Spitzer-Baxter's identity (Proposition \ref{prop Spitzer}), and this is carried out in detail in Section \ref{section GG1} for the correlated GI/G/1 queue and the Sparre-Andersen risk reserve process.
Roughly speaking, all the transforms of the relevant performance measures are obtained by reading the inversion formula in Theorem \ref{generalized Hewitt} from right to left.

Having obtained a Cauchy integral representation for the busy period, one can then evaluate it when the transform of the generic pair $(B,A)$ is a rational function in the argument that corresponds to the service requirement $B$ (Section \ref{section GK1}). 

\vspace{0.3cm}
\section{On Hewitt's inversion formula\label{section inversion}}

The starting point is P. L\'evy's inversion formula which gives a precise form to the well known assertion that a characteristic function uniquely determines a probability measure $\phi$ on the real line:

\[\lim\limits_{T\rightarrow\infty}\frac{1}{2\pi i}\!\int\limits_{-iT}^{iT} \!\left\{\, \int\limits_{-\infty}^\infty\! e^{\xi x} \phi(\de x) \right\} \frac{e^{-\xi a}- e^{-\xi b}}{\xi}\,\de \xi =\frac{1}{2}\phi((a,b)) + \frac{1}{2}\phi((a,b]). \ \]
Thus the value $\phi((a,b]):= \int \chi_{(a,b]}(u)\,\phi(\de u)$ can be recovered, $\chi_{(a,b]}$ being the indicator function of the interval $(a,b]$. Hewitt's formula extends this result to recover directly  functionals of the form

\[\phi(f) := \frac{1}{2}\int [f(u+) + f(u-)] \phi(\de u),\]
for functions $f$ of bounded variation; if $f$ is also continuous, then the integral above becomes $\int f(u)\,\phi(\de u).$


 The inversion formula in the form given by L\'evy is further generalized to higher dimensions and some other topological groups in Hewitt \cite{Hewitt}. It is not however the purpose of the current paper to explore the possibility of a most general form for it. Such an attempt may not even yield satisfactory results, as was already pointed out in \cite{Hewitt}. 

We will generalize the above in Theorem \ref{generalized Hewitt} to probability measures on $\R^2$ related to the random vector $(B,A)$.
The probabilistic structure is never really lost, the conditional distribution of $B$ given $A$ appears throughout the proof of Theorem \ref{generalized Hewitt} disguised as the conditional kernel $q(u,y)$, which is defined below as a "partial" Radon-Nikod\'ym derivative.

\subsection*{Preliminaries}

We say that $f:\R\rightarrow \C$ is of bounded variation if both its real and imaginary parts are of bounded variation; this is the same as $|f|$ being of bounded variation because all norms are equivalent on $\C$. We will also work with (complex-valued) measures on the Borel subsets of $\R^2$, which are of finite total variation. For such a measure $\phi$, we will denote by $|\phi|$ the total variation measure of $\phi$, and with $\|\phi\|$ its total variation. This will suffice for our purposes, but this set-up is fully detailed and generalized in \cite{Hewitt}, and the references therein.

Let $(B,A)$ be a random vector on some probability space, having an arbitrary probability distribution. Denote by $\prob$ the probability measure and by $H$ the joint c.d.f. of $(B,A)$:

\[H(x,y)=\prob (B\leq x,\,A\leq y).\]
\noindent The correlation device for the increment of the random walk is given by this  distribution, not necessarily having a product form, and its Fourier-Stieltjes transform

\[h(s_1,s_2):= \mean e^{-s_1 B-s_2 A}= \int e^{-s_1 x -s_2 y}\,H(\de x,\de y).\]
In general this is convergent only for $\mathcal Re\, s_1 = 0,\,\mathcal Re\,s_2 = 0$, but if $(B,A)$ is supported on the non-negative quadrant in $\R^2$, then $h$ can be continued analytically to $\mathcal Re\, s_1 \geq 0,\,\mathcal Re\,s_2 \geq 0 $.
\noindent The characteristic function of the increment $A-B$ will also be relevant:

\[ h(\xi,-\xi)= \int e^{\xi x} \,\de \prob (A-B\leq x), \,\,\,\mathcal Re\, \xi=0.  \]

Let $\lambda$ be the probability measure associated with the random vector $(B,A)$, and $\nu$ be the marginal measure associated with $B$,

  \[\lambda(\mathcal U\times \mathcal V):= \prob(B\in \mathcal U, A\in \mathcal V),\;\;\;\nu(\mathcal U):=\lambda(\mathcal U\times \R),\;\;\;\mathcal U,\mathcal V\in \borel (\R),\]
 with $\borel(\R)$ the family of Borel sets on the line. The notation $H(\de u, y)$ will be used to suggest that we are integrating w.r.t. the measure $\lambda_y(\mathcal U): = \lambda(\mathcal U\times (-\infty,y]).$
We will work with a version of the {\em conditional} cumulative distribution function (c.d.f.) of $A$ given $B$ and this is made precise below.

  It clearly holds that $\lambda_y(\mathcal U)\leq \nu(\mathcal U)$, in particular $\lambda_y \ll \nu$, so let

  \[q(u,y):=\frac{\de\lambda_y}{\de \nu} (u)\]
  be its Radon-Nikod\'ym derivative. Heuristically, $q(u,y)$ is to be regarded as $q(u,y)=\prob(A\leq y\,|\,B \in \de u)$,
  and we have the disintegration identity

  \[ \int_{\mathcal U} H(\de u,y) = H(\mathcal U,y) =\int_{\mathcal U}q(u,y)\,\nu(\de u),\]
  with any of the terms above meaning $\prob(A\leq y,\;B\in \mathcal U)$.
  We will be working with a \emph{regular version} of $q$, which exists by virtue of the separability of $\R$, see for instance Kallenberg \cite{Kallenberg} Thm. 5.3, p. 84. Furthermore, we have more than just regularity for this kernel, the same result gives that $q(u,y)$ is {\em regularly monotone} as a function in the argument $y$, i.e., $q(u,y)$ is non-decreasing in $y$, for $u$ outside a set of $\nu-$measure zero which does not depend on $y$.

  These considerations are quite intuitive because of the probabilistic nature of the measure associated with $H$. It turns out,  however, that we will have to consider instances of the inversion theorem for the slightly more general case of complex valued functions $H$ which are also of bounded variation, and for this purpose we will show below that $H(\de u,x)$ can be given a meaning in an analogous way.

Confusing $H$ with its associated complex-valued measure, see Hewitt~\cite{Hewitt}, we can reduce it to a monotonically increasing function, by splitting it into real and imaginary parts and using the Jordan decomposition

\[\real H \equiv (\real H)^+ - (\real H)^-,\]
for the signed measure $\real H$, and similarly for $\im H$. Setting $\nu^\pm(\mathcal U):=(\real H)^\pm(\mathcal U\times\R)$, it holds with the similar notation as in the probabilistic case that $(\real H)^\pm_y\ll \nu^\pm$, so we can define again the Radon-Nikod\'ym derivatives

\[q^+_1(u,y):=\frac{\de (\real H)_y^+ }{\de \nu^+}(u),\,\,\,q_1^-(u,y):= \frac{\de (\real H)_y^- }{\de \nu^-}(u),\]
and similarly for $\im H$. Now we can reconstruct $H(\de u,y)$ in an obvious way, using the linearity of the Radon-Nikod\'ym derivative.

Alternatively, we could have used the total variation measures: \[|\real H_y|(\mathcal U)\leq |\real H|(\mathcal U \times\R)\] and thus $\real H_y\ll \nu$, but a simple argument relying on the Hahn-Jordan decomposition shows that this construction yields the same result for $H(\de u,x)$.

Moreover, the monotone regularity property of the probabilistic instance extends to $q(u,y)$ being of bounded variation in $y$ outside a set of $\nu$-measure zero which does not depend on $y$. This property will be useful in the proof of the next result.





\bthm[\textbf{Hewitt inversion extended}]\label{generalized Hewitt} Let $H$ be a totally bounded (complex-valued) measure on $\mathbb R^2$, and let $f:\mathbb R\rightarrow \mathbb C$ be a function of bounded variation which is also absolutely integrable (w.r.t. the Lebesgue measure). Then the following Cauchy principal value can be represented as a Lebesgue-Stieltjes integral:

\[ \lim\limits_{T\rightarrow\infty}\!\frac{1}{2\pi i}\!\int\limits_{-iT}^{iT} \!\left\{ \int\limits_{-\infty}^\infty \int\limits_{-\infty}^\infty \!\! e^{\xi (u-y)}\! H(\de u,\!y) f(y) \de y \!\right\}\!\de \xi  \!=\! \frac{1}{2}\!\int\limits_{-\infty}^\infty \!\left\{f(u+)H(du,\!u+)\!+ \!f(u-)H(du,\!u-)\right\}. \]
\ethm
\noindent $f$ need not be integrable w.r.t. $ H(\de u, u)$. If one of the sides above converges, so does the other one.

\begin{Remark}\label{remark original Hewitt} If $H$ is of the form $H_1H_2$, the double integral inside the Cauchy principal value factorizes into the Fourier-Stieltjes transform of $H_1$ and the Fourier transform of $f(y)H_2(y)$, this function being again absolutely integrable and of bounded variation. Thus the above reduces to the inversion formula in Hewitt~\cite{Hewitt}, Thm. (3.1.1):

\[ \lim\limits_{T\rightarrow\infty}\frac{1}{2\pi i}\int\limits_{-iT}^{iT} \left\{ \int\limits_{-\infty}^\infty  e^{\xi u} H_1(\de u)\right\} \left\{\int\limits_{-\infty}^\infty e^{-\xi y}g(y)\,\de y\right\}\,\de \xi  = \frac{1}{2}\int\limits_{-\infty}^\infty \left[ g(u+)+ g(u-)\right]\,H_1(du),   \]
with $g(y)=f(y)H_2(y).$

\end{Remark}

\proof[Proof of Theorem \ref{generalized Hewitt}]

For fixed $u$, change the variable $x:=u-y$, so that the double integral inside the Cauchy principal value becomes

\[ \int\limits_{-\infty}^\infty \int\limits_{-\infty}^\infty e^{\xi (u-y)} H(\de u,y) f(y)\, \de y = \int\limits_{-\infty}^\infty  \int\limits_{-\infty}^\infty e^{\xi x} H(\de u,u-x) f(u-x)\,\de x . \]

\noindent We can bound

\[ \int \de x\, \left|\int H(\de u, u-x)f(u-x) \right| \leq \int \de x\int |f(u-x)|\, |H|(\de u,\infty) = \|H\|\int |f(v)|\,\de v, \]
hence \noindent $\int\limits  H(\de u,u-x)f(u-x)$ is absolutely integrable in $x$ because $f$ is.
 We can now change the order of integration in the Cauchy principal value, which becomes after integrating over $\xi$:

\begin{equation}\label{principal value intermed} \lim\limits_{T\rightarrow \infty} \iint \frac{1}{2\pi i} \int\limits_{-iT}^{iT} e^{\xi x}\,\de \xi\, H(\de u,u-x) f(x)\, \de x = \lim\limits_{T\rightarrow \infty} \int  \frac{\sin Tx}{\pi x}\de x \!\int\! H(\de u,u-x) f(u-x). \end{equation}


At this point, disintegrate the kernel $H(\de u,u-x) = q(u,u-x)\, \nu(\de u)$, so that, by the Radon-Nikod\'ym Theorem, we can rewrite the right-hand side in (\ref{principal value intermed}) as

\beq  \lim\limits_{T\rightarrow \infty} \int \nu(\de u) \left\{\int \frac{\sin Tx}{\pi x}  q(u,u-x) f(u-x)\, \de x\right\}.\label{R_N Dirichlet}\eeq
The main remark is that  $q(u,u-x)f(u-x)$ is of bounded variation as a function in $x$ for $\nu$-almost all $u$, and regularity is the key, as can be seen from the following:

Let $\{x_i\}_{i\in I}$ be some ordered sequence determined by the edges of an interval partition of $\R$. 
We can write: 

\begin{align*}  \sum_{i\in I} |q(u,u-x_i) f(u-x_i) - &q(u,u-x_{i+1}) f(u-x_{i+1})| \\ &\leq \sum_{i\in I} |q(u,u-x_i)||f(u-x_i) - f(u-x_{i+1})|\\
&+ \sum_{i\in I} |f(u-x_{i+1})| | q(u,u-x_i) - q(u,u-x_{i+1})|,\;\;u\notin \mathcal Q,\\
\end{align*}
where, by virtue of the regularity of $q$, $\mathcal Q$ is a $\nu-$negligible set, outside which $q(u,\cdot)$ is of bounded variation.
Since $\mathcal Q$ does not depend on the choice for the sequence $\{x_i\}$, we can take the supremum over all such sequences, and using that $f$, $q$ are of bounded variation, gives that also $q(u,u-x)f(u-x)$ is of bounded variation in $x$, for $\nu-$almost all $u$.

 We have arrived at the following limit:

\[\lim\limits_{T\rightarrow \infty} \int\limits_{-\infty}^\infty \frac{\sin Tx}{\pi x} \varphi(u,x)\,\de x = \frac{1}{2}[\varphi(u,0+)+ \varphi(u,0-)], \]
 for fixed $u$ and $\varphi(u,x) := q(u,u-x)f(u-x)$. This identity is known as Dirichlet's integral. The integrability condition for the left-hand side that assures the limit exists is that $\varphi(u,\cdot)$ be of bounded variation. As seen from the above, this assumption is only slightly more general than Dirichlet's original condition of monotonicity for $\varphi(u,\cdot)$. See Doetsch~\cite{Doetsch} Ch. 24, or Titchmarsh~\cite{Titchmarsh} 13.2 (the condition is also known as Jordan's test). 

The proof will be complete as soon as we show that the limit in $T$ can be taken inside the $\nu(\de u)$ integral (\ref{R_N Dirichlet}).
Reduce $\varphi(u,\cdot)$ to a positive monotonically decreasing function. Using the second mean value theorem, for any $\beta>0$, we can find $\beta>\alpha>0$ such that

\[ \left|\int\limits_0^\beta \varphi(u,x) \frac{\sin Tx}{x} \de x\right| =\left|\varphi(u,0) \int\limits_0^\alpha \frac{\sin Tx}{x}\,\de x\right| \leq |\varphi (u,0)| \int\limits_0^\pi \frac{\sin x}{x}\,\de x .  \]
 The final upper bound is obtained by a change of variable $x\rightarrow x/T$. 
 The integral along the negative half-axis can be treated in a similar way, and this is sufficient to allow the interchanging of limit and integration in (\ref{R_N Dirichlet}); this together with the fact that Dirichlet's integral identity holds for $\nu$-almost all $u$ completes the proof.

\vspace{.2cm}


\section{The analysis of fluctuations\label{section GG1}}

Having settled the inversion result in the previous section, let us start with the study of the special type of random walk which was briefly described in the introduction.

It is assumed that the random variable $X$ that stands for the generic increment of the random walk $\{S_n\}_{n\geq 0}$

\[S_n = S_0 + \sum_{i=1}^n X_i,\]
 can be written as the difference $B-A$, with $(B,A)$ having some general distribution supported on the non-negative quadrant of $\mathbb R^2$, as in Section \ref{section inversion}.
Moreover, we take the law $\prob$ of $\{S_n\}$
to be the law conditional on $S_0=0$. Also set $b_n= \sum_{i=1}^n B_i$ and $a_n=\sum_{i=1}^n A_i$.


\vspace{.2cm}

In queueing terms, one can think of the pair $(B,A)$ as the service time of a generic customer together with the time until the arrival of the next customer in a GI/G/1 queueing system - which can be always normalized to unit server speed without losing generality.

Let $N$  be a random variable which is distributed as the number of increments before $S_n$ becomes negative for the first time. In the language of the single server queue, $N$ is distributed as the number of customers served during a busy cycle of the server. Then (assuming unit server speed) $b_N$, $a_N$, $-S_N$ stand for the length of the busy period, that of the busy cycle and respectively the idle time of the server.
\newline In terms of insurance and risk theory, this pair can be interpreted as the time elapsed (and hence the premium $B$ gained) since the last claim incurred together with the amount $A$ claimed through an insurance policy.
Then, conditional on starting with 0 initial capital, $b_N$ is the time to ruin, 
$a_N$ is the total amount claimed until ruin (including the claim that causes it) and $-S_N$ is the deficit at ruin.

The representation given below was obtained in Wendel \cite{Wendel_order} \S 4 as an algebraic identity, slightly more general than Spitzer's identity, who originally obtained in \cite{Spitzer57} the representation for the LST of the successive maxima of the partial sums $S_n$ (see also Baxter and Donsker~\cite{Baxter&Donsker} for a similar derivation that holds for L\'evy processes). Theorem 1 and Identity (9) in Kingman \cite{Kingman62a} are much closer to our purposes.
For the sake of completeness, we cite the relevant result in the following proposition


\begin{Proposition}[Spitzer, Wendel, Baxter]\label{prop Spitzer} With the above notations, it holds that

\begin{equation}\label{joint Spitzer} \mean \{z^N e^{-s_1 b_N - s_2 S_N}\} = 1- \exp\left\{ - \sum_{n=1}^\infty  \frac{z^n}{n} \mean [e^{-s_1 b_n-s_2 S_n }1_{\{ S_n< 0\}}]  \right\}, \end{equation}
which is valid for $\mathcal Re\,s_1\geq 0,\, \mathcal Re\, s_2\leq 0,\,|z|\leq 1$.
\end{Proposition}

In Kingman \cite{Kingman62a}, this identity is obtained  by stopping inside the Spitzer-Wendel identity, i.e. in  the three-dimensional space where $(a_n,b_n,S_n)_n$ is evolving, replace the reflecting hyperplane at $x_3=0$ with an absorbing hyperplane. This is formally carried out in \cite{Kingman62a} by replacing the projection operator used by Wendel~\cite{Wendel_order} with an absorption operator. 



\vspace{0.3cm}

We will use Theorem \ref{generalized Hewitt} in conjunction with the version of Spitzer's identity from Proposition \ref{prop Spitzer} to obtain integral representations for the transforms of the busy period, idle period and the number of customers served during a busy period.

Throughout the rest of the paper, we will use the dashed integral sign as a replacement for the cumbersome limit notation
\[  \lim_{T\rightarrow \infty} \int\limits_{-iT}^{iT}  \equiv \mint{-}\limits_{-i\infty}^{i\infty} . \]

\begin{Theorem}\label{Theorem integral representation} We have the following integral representations for $P:=b_N$ and $I:=-S_N$, valid whenever $ \real s> 0,\,|z|< 1$:  

\begin{equation}\label{busy period}
\mean \{z^N e^{-s P}\} = 1- \exp\left\{\frac{1}{2\pi i}\mint{-}\limits_{-i\infty}^{i \infty} \frac{\de \xi}{s-\xi} \left\{\log[1-z h(\xi,0)] - \log[1-zh(\xi,s-\xi)]\right\} \right\},
\end{equation}

\begin{equation}\label{idle period Hewitt}
\mean \{z^N e^{-s I}\} = 1- \exp\left\{\frac{1}{2\pi i}\mint{-}\limits_{-i\infty}^{i \infty} \frac{\de \xi}{s+\xi}\; \log [1-zh(\xi,-\xi)] \right\}.\end{equation}
 Here $\log$ is the principal branch, that has the cut taken along the negative real axis between $0$ and $\infty$,
 so that it admits the  power series representation
 \[ \log \frac{1}{1-r}= \sum_{n=1}^\infty \frac{r^n}{n},\,\,\, |r|< 1.   \]
 \end{Theorem}
\proof
Below we will use an integration by parts argument and for this reason it will be convenient to introduce the function $G(x,y)= \prob(B\leq x, A>y)$. We can write by integrating over the possible values of $b_n$:

\[\prob (b_n < a_n) = \int\limits_{0}^\infty G_{n}(\de x,x)\,\,\mbox{ and }\,\,\prob(b_n\leq a_n) = \int\limits_{0}^\infty  G_{n}(\de x,x-),  \]
where $G_n(x,y) := H^{*n}(x,\infty) - H^{*n}(x,y)$, $H^{*n}$ is the $n$-fold convolution of $H$ with itself and  $G_{n}(\de x,y)$ is the associated integral kernel, as described in Section \ref{section inversion},

\[ G_{n}(\de x,y ) \equiv \prob(b_n \in \de x, \,a_n > y),\;\; G_{n}(\de x,y-) \equiv \prob(b_n \in \de x, \,a_n \geq y) . \]

Let us begin with $(\ref{busy period})$, which means we start with $(\ref{joint Spitzer})$ for $s_2=0$ and $s_1 = s$. The first step is to represent the expected value inside the series in (\ref{joint Spitzer}):

\begin{equation}\label{normalized integral} \frac{1}{2} \mean [e^{-s b_n} (1_{\{S_n\leq 0 \}} + 1_{\{S_n <0\}})] = \int\limits_{0}^\infty  e^{-s x} [\frac{1}{2} G_{n}(\de x,x- )+ \frac{1}{2} G_{n}(\de x,x+)]. \end{equation}
Before we can use the inversion result from Theorem \ref{generalized Hewitt}, define the following $z$-harmonic measure associated to the sequence $G_{n}(\de x,\de y),$ $n\geq 1$:
\
\beq \label{full harmonic measure}H_z^*(\de x,\de y) = \sum_{n=1}^\infty \frac{z^n}{n} G_{n}(\de x,\de y),\eeq
so that, in particular,

\beq \label{log integral kernel}H_z^*(\de x,y)= \sum_{n=1}^\infty \frac{z^n }{n} G_{n}(\de x,y). \eeq
$H_z^*$ is a complex valued measure proper (i.e. it has finite total variation) for $|z|<1$,

\[||H_z^*||\leq \sum_{n=1}\frac{|z|^n}{n}=\log\frac{1}{1-|z|};\]
moreover, its LST equals

\beq\label{logarithmic LST} \int e^{-s_1 x-s_2 y}H_z^*(\de x,\de y) = -\sum_{n=1}^\infty \frac{z^n }{n}h^n(s_1,s_2), \eeq
 since it holds that $G_n(\de x,\de y) = - H^{*n}(\de x,\de y)$, and where the interchanging of the integral with the series is allowed because of absolute integrability:

\[\iint |e^{-s_1 x-s_2 y}|\, |H_z^*|(\de x,\de y) \leq \sum_{n=1}^\infty \frac{|z|^n}{n}\iint H^{*n}(\de x,\de y)= \log\frac{1}{1-|z|}.\]

Now let us use Theorem \ref{generalized Hewitt} for $f(y)= e^{-sy} \chi_{[0,\infty)}(y)$, $\mathcal Re\,s>0$,
which means  we can write via (\ref{normalized integral}) and (\ref{log integral kernel}):

\beq \sum_{n=1}^\infty \frac{z^n}{n}\mean [e^{-sb_n}(1_{\{S_n\leq 0 \}} + 1_{\{S_n <0\}})]=\frac{1}{2\pi i} \mint{-}\limits_{-i \infty}^{i\infty} \left\{ \int\limits_{0}^\infty \int\limits_{0}^\infty e^{\xi x - (s+\xi)y} H_z^*(\de x,y) \, \de y \right\}\,\de \xi.   \label{intermediate Leb-Stieltjes}\eeq

Assume for simplicity  that $\prob(A=B) = 0$ (i.e. $\prob (S_n=0)$ is null for all $n$). This assumption is not essential, see for example the discussion in Cohen \cite{CohenSingleServer}, p. 284. Then the normalized indicator function $\frac{1}{2}( 1_{\{S_n\leq 0 \}}+ 1_{ \{S_n <0\}})$ that appears on the left-hand side of (\ref{normalized integral}) simplifies to $1_{\{S_n<0\}}$.
  Since $H_z^*(\de x,y)$ is of bounded variation in $y$, we can use the integration by parts formula for Lebesgue-Stieltjes integrals, so that (\ref{intermediate Leb-Stieltjes}) becomes

  \[   \frac{1}{2\pi i}  \mint{-}\limits_{-i \infty}^{i\infty} \frac{\de \xi}{s + \xi} \left\{ \int\limits_{0}^\infty  e^{\xi x } \left[H_z^*(\de x,0-) +\int\limits_{y=0}^{\infty}e^{ -(s+\xi)y}H_z^*(\de x,\de y)\right] \right\}  \]

  \begin{equation}\label{b_n integral representation} = \frac{1}{2\pi i}  \mint{-}\limits_{-i \infty}^{i\infty} \frac{\de \xi}{s+\xi}\left[\,\sum_{n=1}^\infty  \frac{z^n}{n}h^n(-\xi,0) - \sum_{n=1}^\infty  \frac{z^n}{n} h^n(-\xi,s+ \xi)\right],   \end{equation}
  where we used (\ref{logarithmic LST}) and the identity $G(x,0-) = \prob (B\leq x)$.
  Changing the variable $\xi\rightarrow -\xi$, the exponent in (\ref{joint Spitzer}) can be rewritten via (\ref{b_n integral representation}), for $\real s>0$, $|z|< 1$: 

\beq\label{convergence issues}-\sum_{n=1}^\infty  \frac{z^n}{n} \mean [e^{-s b_n }1_{\{ S_n< 0\}}] = \frac{1}{2\pi i}   \mint{-}\limits_{-i \infty}^{i\infty} \frac{\de \xi}{s-\xi}\left\{\log[1 - z h(\xi,0)]-\log [1- z h(\xi,s-\xi )] \right\}. \eeq
 Thus (\ref{busy period}) follows from these considerations and Spitzer's identity (\ref{prop Spitzer}).

\vspace{0.2cm}
For the integral representation (\ref{idle period Hewitt}), the extension of Hewitt's formula is not needed. The starting point is (\ref{joint Spitzer}) with $s_1=0$, $s_2=-s$, $\mathcal Re \,s>0$, together with the identity

\[ \mean e^{s S_n}1_{\{S_n<0\}} = \frac{1}{2}\int\limits_{-\infty}^{0+} e^{s x} \de \prob(S_n\leq x)+ \frac{1}{2}\int\limits_{-\infty}^{0-} e^{s x} \de \prob(S_n\leq x) \]
valid because $\prob(S_1=0)=0$. Similarly as (\ref{full harmonic measure}), set $F(x)=\prob(S_1\leq x)$ and introduce

\beq F_z^*(\de x):=\sum_{n=1}^\infty \frac{z^n}{n}F^{*n}(\de x).\label{harmonic measure}\eeq
Use the inversion formula in Remark \ref{remark original Hewitt} with $g(y)=e^{s y}\chi_{(-\infty,0]}(y)$, $\real s>0,$ so that we can write similarly as for (\ref{intermediate Leb-Stieltjes})-(\ref{b_n integral representation}):


\begin{align}\sum_{n=1}^\infty \frac{z^n}{n}\mean e^{s S_n}1_{\{S_n<0\}} = &\frac{1}{2\pi i}  \mint{-}\limits_{-i \infty}^{i\infty} \left\{\int\limits_{-\infty}^{\infty}  e^{\xi x} F_z^*(\de x)\right\} \left\{\int\limits_{-\infty}^{0} e^{(s-\xi)y}\,\de y\right\} \,\de \xi \notag\\ = &\frac{1}{2\pi i}  \mint{-}\limits_{-i \infty}^{i\infty}\frac{\de \xi}{s-\xi} \sum_{n=1}^\infty \frac{z^n}{n} h^n(-\xi,\xi) \notag\\
=&-\frac{1}{2\pi i}  \mint{-}\limits_{-i \infty}^{i\infty}\frac{\de \xi}{s+\xi} \log[1-z h(\xi,-\xi)],\;\;\mathcal Re\, s>0, \label{integral idle period}
 \end{align}
after the change of variable $\xi\rightarrow -\xi$. 

Once (\ref{integral idle period}) is replaced into (\ref{joint Spitzer}), it immediately yields (\ref{idle period Hewitt}). The proof is complete.

  \begin{Remark} If $B$ is independent of $A$, then $h(s_1,s_2)$ is the product of the marginal transforms and  the integral representation  (\ref{busy period}) reduces to that from Kingman~\cite{Kingman62a}, Thm. 4 (see also  Cohen~\cite{CohenSingleServer}, p. 304 for (\ref{idle period Hewitt})).
  \end{Remark}

\begin{Remark}
The integral representations (\ref{busy period}), (\ref{idle period Hewitt}) hold under very general conditions (there are no regularity assumptions for the distribution of $(B,A)$, these can even be discrete random variables, in which case the LSTs become generating functions). The reason is that these representations are given for the interior of their convergence domains ($\real s>0$, $|z|<1$). If we want to take any of the arguments to their respective boundary, we have to require extra conditions to ensure convergence.
For example, when letting $s$ converge towards the imaginary axis, there is a singularity appearing, because the factor $1/(\xi-s)$ gains a simple pole located at $\xi=s$.

It turns out one can give a definite meaning to these integrals, for $\real s=0$, if they are regarded as singular integrals w.r.t. the Cauchy kernel $1/(\xi-s)$. Then one considers the Cauchy principal value obtained by removing a circle of arbitrarily small radius around the singularity $s$ and then taking its radius to 0. Now we are dealing with a double principal value: first one coming from the pole at $s$ of the Riemann integral along the segment $L:=[-iT,iT]$ ($T$ large enough, so that $s\in L$) and the second one obtained by letting $T\rightarrow \infty$. A standard condition (see Gakhov~\cite{Gakhov}, Muskhelishvili~\cite{Mushkelishvili}) that ensures the first principal value converges is that the density functions $\varphi_1(\xi):=\log[1- h(\xi,-\xi)]$, $\varphi_2(s):=\log[1- h(\xi,s-\xi)]$ are H\"older continuous along the imaginary axis, with some positive indices. In the case of (\ref{waiting time Hewitt}), the H\"older continuity of $\varphi_1(\xi)$ is fairly close to Spitzer's~\cite{Spitzer56} integrability condition, which requires (upon taking $s\rightarrow 0$) that $[1-h(\xi,-\xi)]/\xi$ be integrable on a neighbourhood of $s=0$ on the imaginary axis.
\label{remark boundary}
\ermk

%
\vspace{0.2cm}


\subsection*{The number of arrivals during an excursion.}

Further along the lines of Remark \ref{remark boundary}, we will use the doubly dashed integral sign to denote the double Cauchy principal value. 
The choice for the branch of $\log$ is essential for the definiteness of the first principal value. It turns out to be convenient to work with a branch that has the cut between 0 and $\infty$ taken inside the negative half-plane, so we have by definition,

\beq \label{Cauchy pv definition}\frac{1}{2\pi i}\mint{-}\limits_{L}\frac{\varphi(\xi)}{\xi-s}\de\xi  := -\frac{1}{2}\varphi(s) + \frac{\varphi(s)}{2\pi i}[\log(iT-s) - \log(-iT -s)] + \frac{1}{2\pi i} \int\limits_{L} \frac{\varphi(\xi)- \varphi(s)}{\xi-s}\de\xi. \eeq
  The integral on the right is well defined as a Riemann integral as soon as $\varphi(s)$ is H\"older continuous along the line $L$.  By choice of logarithm, the argument of $\log(iT-s)- \log(-iT-s)$ equals $\pi i$ for all $T$. This means that the first two terms cancel in the limit $T\rightarrow \infty$, and the above definition becomes
\beq\frac{1}{2\pi i}\mint{=}\limits_{-i\infty}^{i\infty}\frac{\varphi(\xi)}{\xi-s}\de\xi =\frac{1}{2\pi i} \mint{-}\limits_{-i\infty}^{i\infty} \frac{\varphi(\xi)- \varphi(s)}{\xi-s}\de\xi. \label{Cauchy principal value}\eeq

Formula (\ref{Cauchy pv definition}) differs from the definition given in Mushkelishvili~\cite{Mushkelishvili}, p. 27 or in Gakhov \cite{Gakhov}, p. 16, because therein the cut of the logarithm is taken in the opposite half-plane. To be more precise, in our case, Cauchy's integral representation reads:
\[\mint{-}\limits_{-i\infty}^{i \infty} \frac{\de \xi}{\xi -s} = \left\{\begin{array}{l} 2\pi i,\;s\in L^+, \\ 0,\;s\in L^-, \\ 0,\;s\in L.\end{array}\right. \]

The first two values are the well known Cauchy integral identities; the third one is the Cauchy principal value for this specific choice of $\log$ (using (\ref{Cauchy principal value}) with $\varphi\equiv 1$).

It helps to think about the Riemann sphere as the one point compactification of the complex plane so that the imaginary axis is closed into a large circle on the sphere. Then all of the conventions above are specifying a means of integrating on the large circle of the sphere, the integrands being extended by continuity at infinity; the interior of the imaginary axis is by definition the left hemisphere and the exterior is the right hemisphere.
Define

\[\Phi(s)=\mint{-}\limits_{-i\infty}^{i\infty}\frac{\varphi(\xi)}{\xi-s}\de \xi,\;\; s\in \C,\]
 and the integral is defined in the sense of (\ref{Cauchy principal value}), when $\real s=0$. Let us denote by $\Phi^-(s)$ the limit as $s$ approaches the imaginary axis from its exterior (the right half-plane), and by $\Phi^+(s)$, the limit taken from the interior. Then the Plemelj--Sokhotski formulae (cf. Gakhov \cite{Gakhov}, p. 25) become with the above conventions:

  \beq\Phi^-(s) = \Phi(s),\;\;\; \Phi^+(s) = \Phi(s) + \varphi(s),\;\;\;\real s=0.\label{Plemelj-Sokhotski}\eeq
  Now we can calculate the limit as $s$ approaches the imaginary axis from its exterior in (\ref{busy period}), and  from its interior in (\ref{idle period Hewitt}), using (\ref{Plemelj-Sokhotski}). Still denoting the  interior limit with $\Phi^+(s)$, (\ref{idle period Hewitt}) becomes for $|z|<1$ and as $s$ tends to the imaginary axis:

 \beq\label{phi plus}\Phi^+(s) = 1 -[1-z h(- s,s)]\exp\left\{  \frac{1}{2\pi i}\mint{=}\limits_{-i\infty}^{i \infty} \frac{\de \xi}{\xi+s}\; \log [1-z h(\xi,-\xi)] \right\}.\eeq
Before we give the limiting values of (\ref{busy period}), we point out how one can simplify it. Consider the first integral in the exponent of (\ref{busy period}):

\beq\label{first term}\Psi(\eta):=-\frac{1}{2\pi i}\mint{-}\limits_{-i\infty}^{i \infty} \frac{\de \xi}{\xi-\eta} \log[1-z h(\xi,0)],\;\;\real \eta>0.\eeq
 The key remark is that $h(\xi,0)=\mean e^{-\xi B}$ is analytic for $\real \xi>0$ and this implies that the density $\log[1-zh(\xi,0)]$ is again analytic for $\real \xi>0$, because of the choice of $\log$. 
Since the integrand has a simple pole in the positive half-plane located at $\eta=\xi$, it follows at once from Cauchy's Theorem applied to the imaginary axis that $\Psi(\eta)$ equals

\[\Psi(\eta)=\log[1-z h(\eta,0)],\]
for any $\real \eta>0$. This can be seen by closing the segment $[-iT,iT]$ with the half circle spanning between its endpoints inside the positive half-plane. The contour integral thus obtained equals $\log[1-zh(s,0)]$ for $T$ large enough, so that the pole $\eta=\xi$ lies inside the contour (remark that because of the conventions on the interior of the imaginary axis, this contour is traversed in the clockwise direction). Finally, the contribution along the half-circle tends to 0 as $T\rightarrow \infty$ because the integrand behaves as $o(|\xi|^{-1})$ along the half-circle.

 Having settled the first term in the exponent of $(\ref{busy period})$, we can use the Plemelj--Sokhotski's formula (the continuity of the exterior limit from (\ref{Plemelj-Sokhotski})) for the other term, to obtain the identity

\beq \Phi^-(s) = 1- \exp\left\{  \log[1-z h(s,0)] + \frac{1}{2\pi i}\mint{=}\limits_{-i\infty}^{i \infty} \frac{\de \xi}{\xi-s}\; \log [1-zh(\xi,s-\xi)] \right\},\; \real s=0,\label{phi minus 0}\eeq
hence we can rewrite (\ref{phi minus 0}) as

\beq\label{phi minus}  \Phi^-(s) = 1- [1-zh(s,0)]\exp \left\{ \frac{1}{2\pi i}\mint{=}\limits_{-i\infty}^{i \infty} \frac{\de \xi}{\xi-s}\; \log [1-zh(\xi,s-\xi)] \right\},\;\; \mathcal Re\, s=0.\eeq


 In particular, for $s=0$, the limits (\ref{phi plus}) and (\ref{phi minus}) agree and these must then coincide with the generating function of $N$. The following has been proven

\bprop \label{Prop N}With the above notations and conventions, it holds for $|z|<1$:

\beq\mean z^{N} = 1 - (1-z)\exp\left\{  \frac{1}{2\pi i}\mint{-}\limits_{-i\infty}^{i \infty} \frac{\de \xi}{\xi}\; \log \frac{1-zh(\xi,-\xi)}{1-z} \right\}.\label{gen funct N}\eeq
\eprop

One can determine from Spitzer-Baxter's identity (see Spitzer \cite{Spitzer56}, Thm. 3.1) the transform of the transient waiting time in a similar way as above. The Spitzer-Baxter identity reads for $M_n=\max\{S_0,S_1,...,S_n\}$:

\beq \sum_{n=0}^\infty z^n \mean e^{-s M_n} = \exp\left\{ \sum_{n=1}^\infty \frac{z^n}{n} \mean e^{-s S_n^+}\right\}, \label{Spitzer-Baxter}\eeq
 with $S_n^+=\max(S_n,0)$. Keeping the assumption that $\prob(S_1=0)=0$, we can write 

\beq \label{intermediate Spitzer-Baxter}\sum_{n=1}^\infty\mean e^{-sS_n^+} = \int e^{-sx}\chi_{[0,\infty)}(x) \de F^*_z(x) + \sum_{n=1}^\infty \frac{z^n}{n}\int \chi_{(-\infty,0)}(x)\de\prob(S_n\leq x).\eeq
for $F^*_z$ defined in (\ref{harmonic measure}).
It is worth to point out that in this case $S^+_n$ has the mass of $S_n$ distributed along $(-\infty,0]$ swept into the origin as an atom. 
The indicator function appearing in the last term in (\ref{intermediate Spitzer-Baxter}) is not absolutely integrable w.r.t. the Lebesgue measure, therefore  introduce $g_\epsilon(x) = e^{\epsilon x}\chi_{(-\infty,0)}(x)$, $\epsilon>0$, in order to be able to use Hewitt's inversion formula. This same perturbation of the indicator function was used by Spitzer \cite{Spitzer57}. We have similarly as (\ref{integral idle period}):

\[\sum_{n=1}^\infty \frac{z^n}{n}\mean e^{\epsilon S_n}1_{\{S_n<0\}} = -\frac{1}{2\pi i}  \mint{-}\limits_{-i \infty}^{i\infty}\frac{\de \xi}{\xi+\epsilon} \log[1-z h(\xi,-\xi)].\]
We can take $\epsilon\rightarrow 0$, by using the Plemelj-Sokhotski identity (\ref{Plemelj-Sokhotski}) for the right-hand side and dominated convergence on the left ($|z|<1$):

\beq \sum_{n=1}^\infty \frac{z^n}{n}\prob(S_n<0) =  -\frac{1}{2\pi i}\mint{-}\frac{\de \xi}{\xi}\; \log \frac{1-zh(\xi,-\xi)}{1-z}. \label{term two Spitzer-Baxter}\eeq
Hewitt's inversion formula is directly applicable for the first term in (\ref{intermediate Spitzer-Baxter}): as in the proof of (\ref{idle period Hewitt}), use Remark \ref{remark original Hewitt}  for $g(y)=e^{-s y} \chi_{[0,\infty)}(y),$ $\real s> 0$:

\[ \sum_{n=1}^\infty \frac{z^n}{n}\mean e^{-sS_n^+} = \frac{1}{2\pi i} \mint{-}\limits_{-i\infty}^{i \infty}\left\{\sum_{n=1}^\infty \frac{z^n}{n}\mean e^{\xi S_n} \int\limits_{0}^\infty e^{-(s+\xi)y}\de y \right\}\de \xi = \frac{1}{2\pi i} \mint{-}\limits_{-i\infty}^{i \infty} \frac{\de\xi}{s-\xi} \sum_{n=1}^\infty \frac{z^n}{n}h^n(\xi,-\xi),\]
after the change of variable $\xi\rightarrow -\xi$. This together with (\ref{Spitzer-Baxter}) yields

\beq \label{waiting time Hewitt} \sum_{n=0}^\infty z^n \mean e^{-s M_n}\! =\!\frac{1}{1-z}\exp\!\left\{ \!\frac{1}{2\pi i}\!\mint{-}\limits_{-i\infty}^{i \infty}\! \frac{\de \xi}{\xi-s}\! \log [1\!-\!z h(\xi,-\xi)]\!- \!\frac{1}{2\pi i}\mint{-}\limits_{-i\infty}^{i \infty}\! \frac{\de \xi}{\xi} \log \frac{1\!-\!z h(\xi,-\xi)}{1-z}\!\right\}\eeq
with $\real s>0,\;|z|<1$ (cf. Cohen \cite{CohenSingleServer}, (5.29) p. 276).
This formula is related to that obtained by Spitzer \cite{Spitzer57} concerning the Wiener-Hopf equation which has a probability density as the kernel.

\section{Examples\label{section GK1}}


In this section we evaluate the integral representation (\ref{busy period}), under the assumption that the transform of the generic pair $(B,A)$ is a rational function in the argument that corresponds to the service requirement $B$. The analysis is similar to the one carried out in Cohen~\cite{CohenSingleServer}, Ch. II.5.



Assume that for all $s_2$, the joint LST $h(s_1,s_2)$ is a rational function in the argument $s_1$, which can be represented as

\beq h(s_1,s_2) = \frac{h_1(s_1,s_2)}{h_2(s_1,s_2)}, \label{rational definition} \eeq
where $h_i(\cdot,s_2)$ are polynomial functions. Moreover we will assume that for $\real s\geq 0$, $h(\xi,s-\xi)$ has a finite number of poles in the negative half-plane as a function in the argument $\xi$ ($h(\xi,s-\xi)$ is already meromorphic in this region, because of the previous assumption). 
This is not essential, but an algorithmically friendly assumption which will give a representation for the busy period transform in terms of a finite number of factors.

We may still assume without losing generality that $\prob(B-A=0)=0$, which implies $h(s_1,s_2)\rightarrow 0$, as $s_1\rightarrow \infty$, $\real s_1>0$, and the convergence is uniform in $s_2$, for $\real s_2\geq 0$.
In particular, we have for any $\real s_2\geq 0$, $\deg h_1(\cdot,s_2)<\deg h_2(\cdot,s_2)$. 

Before we proceed with the analysis, let us point out some ways of creating correlation between the inter-arrivals and the corresponding service times.

\vspace{.2cm}
\textbf{Example 1} (Threshold dependence) This is one of the simplest ways of making $B$ depend on the size of $A$: for a fixed threshold $l>0$, $B\sim B_1$ on the event $A\leq l$ and $B\sim B_2$ otherwise; with $B_i$ independent of $A$ and having rational transforms $f_i(s_1)$, $i=1,2$; thus

\[h(s_1,s_2) = f_1(s_1)a_1(s_2) + f_2(s_1) a_2(s_2),  \]
\[a_1(s_2) = \int\limits_{y=0}^l e^{-s_2 y} \,\de\prob(A\leq y), \;\;\;\;\;\;a_2(s_2) = \int\limits_{y=l+}^\infty e^{-s_2 y} \,\de\prob(A\leq y), \]
so that $a_1(s_1)$ is an entire function and $a_2(s_2)$ is analytic and bounded for $\real s_2>0$. This construction can be naturally extended to $k$ thresholds, giving

\[h(s_1,s_2)= \sum_{i=1}^k f_i(s_1) a_i(s_2),\]
with $a_i(s_2)$ entire functions, $i<k$, and $a_k(s_2)$ is analytic for $\real s_2>0$.

\vspace{.2cm}
\textbf{Example 2} (Markov Modulation) Let $(X_n)$ be a finite state Markov chain which has an absorbing state and denote with $\kappa$ the number of jumps until absorbtion. Define   

\[(A,B) = \sum_{i=1}^\kappa (A_i,B_i),\]
where $(A_i,B_i)$ are i.i.d. vectors. The component $A_1$ is allowed to be generally distributed with $g_0(s_2)=\mean e^{-s_2 A_1}$,  and $B_1$ has a rational transform of the form $f_1(s_1)/f_2(s_1)$. 

If we denote by $\boldsymbol\alpha$ the initial distribution of $(X_n)$, by $T$ the transient component of its transition matrix, and by $\boldsymbol t$ the vector of exit probabilities, then by conditioning on $\kappa$, the transform of $(A,B)$ is (see for instance \cite{paper1}):

\[ h(s_1,s_2) = \boldsymbol\alpha^t\left[\frac{f_2(s_1)}{f_1(s_1)g_0(s_2)}I-T\right]^{-1}\boldsymbol t.  \]
Both these examples are of the form assumed by (\ref{rational definition}).

\vspace{.2cm}
Remark that $h_2(\cdot,s_2)$ can only have zeroes with negative real part, because of the regularity domain of $h(\cdot,s_2)$.
With these assumptions, the exponent in (\ref{busy period}) becomes

\beq \log[1-\mean z^Ne^{-sP}]=\frac{1}{2\pi i}\mint{-}\limits_{-i\infty}^{i \infty} \frac{\de \xi}{s-\xi}\log[1 - z h(\xi,0)] - \frac{1}{2\pi i}\mint{-}\limits_{-i\infty}^{i \infty} \frac{\de \xi}{s-\xi}\log[1-z h(\xi,s-\xi)] \label{rational rep}. \eeq

For the principal branch of the logarithm which has the cut taken along the negative real axis, the single valued functions $\log[1-z h(\xi,0)]$ and $\log[1-z h(\xi,s-\xi)]$  are {\em holomorphic}, for $\xi$ lying in a neighbourhood of infinity, $\real \xi<0$. The reason is that for such values of $\xi$, $|zh(\xi,0)|<1$, $|zh(\xi,s-\xi)|<1$, and then a simple geometric argument shows that both  $1-z h(\xi,0)$ and $1-zh(\xi,s-\xi)$ lie in the positive half-plane.
With this choice of the cut, the evaluation of the integrals (\ref{rational rep}) becomes an application of the theorem of residues. 
Before we can evaluate (\ref{rational rep}), the zeroes and poles of the arguments of the logarithm must be localized. 
The following lemma will also be useful later on.

\blem \label{lemma Rouche Hewitt}The functions $h_2(\xi,s-\xi)$ and $h_2(\xi,s-\xi) - z h_1(\xi,s-\xi)$ have the same number $n\equiv n(s,z)$ of zeroes in the negative half of the complex $\xi$-plane, either when $|z|<1$, $\real s\geq 0$, or $|z|\leq 1$, $\real s>0$.

Assuming $\real s\geq 0$ and $\mean B,\, \mean A<\infty$, then under the extra ergodicity condition $\mean B<\mean A$,  the functions $h_2(\xi,s-\xi)$ and $h_2(\xi,s-\xi) -  h_1(\xi,s-\xi)$ have the same number $m\equiv m(s)$ of zeroes with negative real part.
\elem

\proof  The statements will follow from Rouch\'e's theorem (cf. Titchmarsh \cite{Titchmarsh}, p. 116) as soon as we show that these functions are analytic in the interior of some suitably chosen contours and on their boundary it holds that

\beq | h_2(\xi,s-\xi)| >|z h_1(\xi,s-\xi)|,\;\;\;\; |h_2(\xi,-\xi)| > |h_1(\xi,-\xi)|. \label{contour inequalities}\eeq

Fix $R>0$ and consider the contour $\mathcal C$ consisting of the segment of the imaginary axis between $-iR$ and $iR$ together with the semicircle with radius $R$ that spans in the negative half-plane.
For the segment of the imaginary axis, we have the following bounds

\beq  |z h(\xi,s-\xi)|= |z|| \mean e^{-\xi B - (s-\xi)A}|\leq  |z|\,\mean |e^{-\xi B - (s-\xi) A}| \leq |z|\,\mean e^{- \real s A}, \;\;\; \mathcal Re\, \xi = 0.  \label{imaginary bound} \eeq
For the bound on the half-circle, consider the following representation for $h(s_1,s_2)$:

\[h(s_1,s_2) = \frac{a_1(s_2) s_1^{n-1} + a_2(s_2)s_1^{n-2}+\ldots+a_n(s_2)}{b_1(s_2)s_1^{n}+ b_2(s_2)s_1^{n-1} \ldots+b_{n+1}(s_2)},\]
where $n\equiv n(s_2)$ (remark that for the examples presented above, the denominator does not depend on $s_2$, hence neither does the degree $n$). The functions $a_i(s_2)$ can be taken to be bounded for $\real s_2 \geq 0$, because it holds that $|h(1,s_2)|\rightarrow 0$ as $s_2\rightarrow \infty$. Moreover, we can assume $b_1(s_2)\equiv 1$, for $\real s_2\geq 0$, after normalizing the fraction. For fixed $s_2$, let $\xi_i(s_2)$ be the zeroes (all with negative real part) of $h_2(s_1,s_2)$; when bounding the above representation of $h(s_1,s_2)$, use the triangle inequality for the numerator and use the inequality  $|z_1-z_2|\geq ||z_1|-|z_2||$ for the denominator:

\[|h_2(s_1,s_2)|= \left|\prod_{i=1}^n(s_1 - \xi_i(s_2))\right|\geq \left|\prod_{i=1}^n(|s_1|-|\xi_i(s_2)|)\right|,\]
the right-hand side being a polynomial function in $|s_1|$ of the same degree as $h_2(s_1,s_2).$ Thus we have the upper bound

\[ |h(\xi,s-\xi)|\leq \frac{|a_1(s-\xi)|\,|\xi|^{n-1}+ |a_2(s-\xi)|\,|\xi|^{n-2}+\ldots}{|\prod_{i=1}^n(|\xi|-|\xi_i(s-\xi)|)|}.\]

Then it follows from the facts that $\deg h_1(\xi,s-\xi)<\deg h_2(\xi,s-\xi)$ and that the $a_i$ are bounded, that

\beq \label{semi-circle bound} |h(\xi,s-\xi)| = o(|\xi|^{-1}), \;\; |\xi|=R \rightarrow \infty,\, \real s \geq 0,\eeq
for $\xi$ running along the half-circle that closes the contour $\mathcal C$.

 From the bounds (\ref{imaginary bound}) and (\ref{semi-circle bound}), it follows that $|z h_1(\xi,s-\xi)|<|h_2(\xi,s-\xi)|$ when $\xi$ is on $\mathcal C$, for $R$ large enough, either if $|z|<1$, $\real s\geq 0$, or if $|z|\leq 1$, $\real s >0$. This yields the first part of the lemma, via Rouch\'e's theorem.  


For the second part, consider the contour $\mathcal C_\epsilon$ made up from the segment that runs in parallel to the imaginary axis and lying to its left at distance $\epsilon$, together with the arc of the circle with radius $R$ spanning in the negative half-plane between the edges of this segment.

It is essential that $ h(\xi,-\xi)$ is meromorphic in $\real \xi<0$ (see the discussion below (\ref{rational definition})). Since it has isolated poles, we can find $\epsilon>0$, such that $h(\xi,-\xi)$ is holomorphic in the thin strip $-2\epsilon < \real \xi<0$.
Then the left derivative of the function $h(\xi,-\xi)$ exists at 0 and we have by hypothesis,

\[ \lim_{\stackrel{\xi\rightarrow 0}{\real \xi<0}}\frac{\de }{\de \xi} h(\xi,-\xi) = \mean A - \mean B >0, \]
in particular, $h(\real \xi,-\real\xi)<h(0,0)=1$.
We can now bound for $\real s\geq 0$ and $\xi$ lying on the segment of $\mathcal C_\epsilon$:

 \[|h(\xi,s-\xi)| \leq \mean |e^{-\xi B -s A + \xi A}| \leq h(\real \xi,-\real\xi)<1.\]
The bound for $\xi$ lying on the arc component of $\mathcal C_\epsilon$ follows in the same way as (\ref{semi-circle bound}). By virtue of Rouch\'e's theorem, the proof is complete.



\brmk  \label{remark Rouche Hewitt}It follows in a similar way as for the first part of Lemma \ref{lemma Rouche Hewitt} that the polynomials $h_2(\xi,0)$ and $h_2(\xi,0) -z h_1(\xi,0)$ have the same number of zeroes with negative real part, $|z|<1$. But since $h_2(\xi,0)$ has only such zeroes and
$\deg h_2(\cdot,0) >\deg h_1(\cdot,0)$, the same holds for $h_2(\xi,0) - zh_1(\xi,0)$.
\ermk

The idea for evaluating (\ref{rational rep}) is to use the theorem of residues for the contour integrals along $\mathcal C_\epsilon$ while arguing that the contributions from the integrals along the half-circle vanish as the radius $R\rightarrow \infty$.
Focus on the contour integral of the second term in (\ref{rational rep}) taken along the semi-circle component, say $\mathcal S_\epsilon$, of $\mathcal C_\epsilon$. For $R$ large enough $\mathcal S_\epsilon$ will be contained in the interior of a domain where $h(\xi,s-\xi)$ is holomorphic, and in addition, $|zh(\xi,s-\xi)|\leq 1$, hence the position vector $1-z h(\xi,s-\xi)$ has positive real part when the argument $\xi$ runs along $\mathcal S_\epsilon$. This means $\log[1- z h(\xi,s-\xi)]$ is holomorphic in a neighbourhood around the arc $\mathcal S_\epsilon$. In conclusion, we can integrate by parts:

\begin{align*} \int_{\mathcal S_\epsilon}  \frac{\de \xi}{s-\xi}\log[1-z h(\xi,s-\xi)] &= -\log(s-\xi )\log[1-z h(\xi,s-\xi)]\biggr|_{-\epsilon-iR}^{-\epsilon +iR} \\+ &\int_{\mathcal S_\epsilon}  \log(s-\xi)\frac{\de}{\de\xi}\log[1-z h(\xi,s-\xi)]\de\xi. \end{align*}
For large $R$, $|h(\xi,s-\xi)|\rightarrow 0$, which means $|\log[1-z h(\xi,s-\xi)]|\sim |z h(\xi,s-\xi)|$, so the first term on the right behaves in absolute value as

\[ \sim |\log(\xi-s)|\, |z h(\xi,s-\xi)| \sim |z\log(\xi-s)|\,|\xi|^{-1}, \]
and the integrand on the left-hand side behaves as $|\xi|^{-2}$. Thus we have

\beq \int_{\mathcal S_\epsilon} \log(\xi-s)\frac{\de}{\de\xi}\log[1-z h(\xi,s-\xi)]\de\xi \rightarrow 0,\;\;\;\;|\xi|\rightarrow \infty,\;\real \xi<0.\label{circular contribution}\eeq

%

Now we are ready to calculate the contour integrals (\ref{rational rep}).
Remark that the first term is the same as (\ref{first term}) and thus equals $\log[1-z h(s,0)]$, as was found in Section \ref{section GG1}.
For the second integral in (\ref{rational rep}), fix $\real s>0$ and consider the integral taken along the contour $\mathcal C_\epsilon$ described in the proof of Lemma \ref{lemma Rouche Hewitt}. $\epsilon$ is taken sufficiently small such that no poles of the integrands are lying between the segment and the imaginary axis, irrespective of $R$ (this can be found since there are finitely many poles in the negative half-plane). The integral can be approximated from the interior of the negative half-plane using the contours $\mathcal C_\epsilon$, for arbitrarily large $R$ and small $\epsilon$.
Splitting the integral along $\mathcal C_\epsilon$ based on the factors inside the logarithm and integrating by parts in each term, the expression in (\ref{rational rep}) becomes

\[ \frac{1}{2\pi i}\int_{\mathcal C_\epsilon}\log(s -\xi)\,\frac{\frac{\de}{\de\xi} [h_2(\xi,0) - z h_1(\xi,0)]}{h_2(\xi,0) - z h_1(\xi,0)}\de \xi - \frac{1}{2\pi i}\int_{\mathcal C_\epsilon} \log(s -\xi) \frac{\frac{\de}{\de\xi}\, h_2(\xi,0)}{h_2(\xi,0)}\de\xi \]

\[+\frac{1}{2\pi i}\int_{\mathcal C_\epsilon} \log(s-\xi) \frac{\frac{\de}{\de\xi}\, h_2(\xi,s-\xi)}{h_2(\xi,s-\xi)}\de\xi  - \frac{1}{2\pi i}\int_{\mathcal C_\epsilon}\log(s -\xi)\,\frac{\frac{\de}{\de\xi} [h_2(\xi,s-\xi) - z h_1(\xi,s-\xi)]}{h_2(\xi,s-\xi) - z h_1(\xi,s-\xi)}\de \xi.    \]

By (\ref{circular contribution}), the total contribution from the integral along the semi-circle $\mathcal S_\epsilon$ vanishes as $R\rightarrow \infty$.
 Moreover, the branch of $\log$ was chosen such that the factors $\log(s-\xi)$ are analytic for $\real \xi <0$. Then the integrands have simple poles located at the zeroes of their denominators in the negative half of the complex plane. 
 Denote by $\xi_i(s),$ $\xi_i(z,s)$, $i=1,...,n,$ the zeroes with negative real part of $h_2(\xi,s-\xi)$, respectively $h_2(\xi,s-\xi) - z h_1(\xi,s-\xi)$, as functions of the variable $\xi$ (see Lemma \ref{lemma Rouche Hewitt}), and with $\eta_j,$ $\eta_j(z)$, $j=1,...,m$, the zeroes (all having negative real part)  of $h_2(\xi,0)$, respectively $h_2(\xi,0)-z h_1(\xi,0)$.
Then the integral in (\ref{rational rep}) equals

\[\log\frac{  \prod_{j=1}^m[s-\eta_j(z)] \prod_{i=1}^n[s-\xi_i(s)]}{ \prod_{j=1}^m[s-\eta_j] \prod_{i=1}^n[s-\xi_i(z,s)]} = \log\frac{[h_2(s,0)-z h_1(s,0)]\prod_{i=1}^n[s-\xi_i(s)]}{h_2(s,0)\prod_{i=1}^n[s-\xi_i(z,s)]},\]
after letting $R\rightarrow \infty$ and $\epsilon\rightarrow 0.$ Remark also that the degree $n$ is a (piecewise constant) function of the argument $s$ (Lemma \ref{lemma Rouche Hewitt}). In conclusion, we have

 \beq \mean \{z^N e^{-s P}\} = 1- [1-z h(s,0)]\frac{\prod_{i=1}^n[s-\xi_i(s)] }{\prod_{i=1}^n[s-\xi_i(z,s)]}.\label{joint busy period}\eeq

The calculations that led to (\ref{joint busy period}) can be repeated for $z=1$ and $\real s\geq 0$. The contour of integration is the same as $\mathcal C_\epsilon$, and the second part of Lemma \ref{lemma Rouche Hewitt} must be used to conclude about the number of zeroes under the condition $\mean A >\mean B$. This condition is necessary for stability ($P$ has a proper probability distribution, which must be verified when taking $z=1$, $s=0$).
The conclusion is that we are allowed to formally replace $z=1$ in (\ref{joint busy period}), which becomes

\[\mean e^{-s P} = 1- [1- h(s,0)]\frac{ \prod_{i=1}^n[s-\xi_i(s)] }{\prod_{i=1}^n[s-\xi_i(1,s)] },\]
with the remark that $P$ has indeed a proper probability distribution.

Similarly, the exponent in (\ref{gen funct N}), which also appears in (\ref{waiting time Hewitt}), can be rewritten as

\[\frac{1}{2\pi i}\mint{-} \frac{\de \xi}{\xi} \log \frac{1-z h(\xi,-\xi)}{1-z} = \frac{1}{2\pi i}\mint{-} \log(-\xi)\frac{\frac{\de}{\de \xi}[1-z h(\xi,-\xi)]}{1-z h(\xi,-\xi)}\de \xi, \]
whereas the other integral that appears in (\ref{waiting time Hewitt}) becomes
\[\frac{1}{2\pi i}\mint{-}\limits_{-i\infty}^{i \infty}\frac{\de \xi}{\xi-s} \log [1-z h(\xi,-\xi)] = \frac{1}{2\pi i}\mint{-} \log(s-\xi)\frac{\frac{\de}{\de \xi}[1-z h(\xi,-\xi)]}{1-z h(\xi,-\xi)}\de \xi. \]
Using the assumption (\ref{rational definition}), it is easy to see that the denominator of both integrands above is of the form $h_2(\xi,-\xi)[h_2(\xi,-\xi)-z h_1(\xi,-\xi)].$
A similar analysis as for (\ref{joint busy  period}) yields the generating function of the sequence $\{\mean e^{-sM_n}\}_n$ and the generating function of $N$ from Proposition \ref{Prop N}:

\beq (1-z)\sum_n z^n\, \mean e^{-sM_n} = 
\frac{\prod_{j=1}^n(s-\xi_j(0))}{\prod_{j=1}^n(s-\xi_j(z,0))}
\frac{ \prod_{j=1}^n\xi_j(z,0)}{\prod_{j=1}^n \xi_j(0)}, \label{transient waiting time}
\eeq

\[\mean z^N = 1- (1-z)\frac{\prod_{j=1}^n \xi_j(0)}{\prod_{j=1}^n \xi_j(z,0)},\;\;|z|<1,\;\real s>0.\]
$\xi_j(s)$, $\xi_j(z,s)$ are defined as in (\ref{joint busy  period}) and $n\equiv n(s,z)$ is the number of zeroes with negative real part, as given by Lemma \ref{lemma Rouche Hewitt}. Under the stability condition $\mean B < \mean A$, $\lim_{n\rightarrow \infty} M_n = M$ exists in distribution and we can take $z\rightarrow 1$ in (\ref{transient waiting time}) by virtue of Abel's theorem (Titchmarsh \cite{Titchmarsh}, 1.22):


\[\mean e^{-s M} = \frac{\prod_{j=1}^n(s-\xi_j(0))}{\prod_{j=1}^n(s-\xi_j(1,0))}
\frac{ \prod_{j=1}^n \xi_j(1,0)}{\prod_{j=1}^n \xi_j(0)}. \]

Finally, it is possible to derive the transform of the idle period using a similar contour integral; and the methods of this section equally apply to the symmetric case when the transform $h(s_1,s_2)$ is a rational function in the argument $s_2$, for each fixed $s_1$, $\real s_1\geq0$. 
The analysis relies in this case on localizing the poles and zeroes lying inside the positive-half plane.

%
%
%
%

\vspace{.2cm}
\textbf{Concluding remarks}
The methods used in Sections \ref{section GG1} and \ref{section GK1} have been first applied in Queueing theory by Pollaczeck \cite{Pollaczek} and many others have followed (see also the survey of Tak\'acs \cite{Takacs76}).
With respect to calculating the quantities related to the fluctuations of random walks, from a general perspective, the analysis presented in this paper shows the effectiveness of the Radon measure and its associated integral (what we called the Lebesgue-Stieltjes integral) when combined with Spitzer's identity (\ref{joint Spitzer}), via the harmonic measures (\ref{full harmonic measure}) and (\ref{harmonic measure}).



\vspace{.2cm}
The following functions:

\[\psi_z^+(s):=(1-z)\sum z^n \mean e^{-sM_n},\;\;\;\;\real s\geq 0, \;|z|<1,\] and
\[\psi_z^-(s):= 1-\mean z^N e^{s I}, \;\;\;\;\real s\leq 0, \;|z|<1,\]
make up the solution of the homogeneous Hilbert problem associated with the Wiener-Hopf equation. To be more precise, $\psi_z^+(s)$ and $\psi_z^-(s)$ satisfy the boundary relation:
\[\psi_z^-(s)\,\psi_z^+(s) = 1- z h(s,-s),\;\;\real s=0,\,|z|<1,\]
and are unique with the additional property that $\psi_z^+(s)$ is analytic in $\real s>0$ and $\psi_z^-(s)$ is analytic in $\real s<0$, both continuous up to the boundary (the imaginary axis). The uniqueness holds because the kernel $1-z h(s,-s)$ has index zero along the imaginary axis (i.e., its logarithm is single-valued along the imaginary axis, see Gakhov \cite{Gakhov}, Ch. II), a fact which was used throughout Sections \ref{section GG1} and \ref{section GK1}.
It is easy to check that the boundary relation holds, if we use the current version of the Plemelj-Sokhotski identities (\ref{Plemelj-Sokhotski}), and take $s$ to the imaginary axis both in (\ref{idle period Hewitt}) and (\ref{waiting time Hewitt}). 

\vspace{.2cm}
There is a large amount of literature on Wiener-Hopf equations both in Analysis and Probability. A reference in analysis is Krein's manuscript \cite{Krein62}; for probabilistic background and applications, see Asmussen \cite{Asmussen2003}. It seems  Rapoport \cite{Rapoport} was the first to observe the connection between the Wiener-Hopf equation and the theory of Riemann/Hilbert boundary-value problems.

\section*{Acknowledgements}

I am indebted to Onno Boxma for many helpful discussions and suggestions, and in particular for pointing out the usefulness of the Plemelj-Sokhotski identities.

\end{document}